\documentclass{amsart}
\usepackage{amssymb}
\usepackage{eucal}
\usepackage{amsmath}
\usepackage{amscd}

\begin{document}

\newcommand{\delete}[1]{}

\newtheorem{theorem}{Theorem}[section]
\newtheorem{prop}[theorem]{Proposition}
\newtheorem{defn}[theorem]{Definition}
\newtheorem{lemma}[theorem]{Lemma}
\newtheorem{coro}[theorem]{Corollary}
\newtheorem{claim}{Claim}[section]
\newtheorem{remark}{Remark}[section]

\newcommand{\nc}{\newcommand}

\renewcommand{\Bbb}{\mathbb}
\renewcommand{\frak}{\mathfrak}

\newcommand{\efootnote}[1]{}

\renewcommand{\textbf}[1]{}

\nc{\mlabel}[1]{\label{#1}}  

\nc{\bin}[2]{ (_{\stackrel{\scs{#1}}{\scs{#2}}})}  
\nc{\binc}[2]{ \left (\! \begin{array}{c} \scs{#1}\\
    \scs{#2} \end{array}\! \right )}  
\nc{\bincc}[2]{  \left ( {\scs{#1} \atop
    \vspace{-1cm}\scs{#2}} \right )}  
\nc{\disp}[1]{\displaystyle{#1}}
\nc{\scs}[1]{\scriptstyle{#1}}
\nc{\proofbegin}{\begin{proof}}
\nc{\proofend}{\end{proof}}
\nc{\sha}{{\mbox{\cyr X}}}  
\nc{\shpr}{\diamond}    
\nc{\vep}{\varepsilon}
\nc{\labs}{\mid\!}
\nc{\rabs}{\!\mid}

\nc{\BA}{{\Bbb A}}
\nc{\CC}{{\Bbb C}}
\nc{\DD}{{\Bbb D}}
\nc{\EE}{{\Bbb E}}
\nc{\FF}{{\Bbb F}}
\nc{\GG}{{\Bbb G}}
\nc{\HH}{{\Bbb H}}
\nc{\LL}{{\Bbb L}}
\nc{\NN}{{\Bbb N}}
\nc{\QQ}{{\Bbb Q}}
\nc{\RR}{{\Bbb R}}
\nc{\TT}{{\Bbb T}}
\nc{\VV}{{\Bbb V}}
\nc{\ZZ}{{\Bbb Z}}


\nc{\cala}{{\cal A}}
\nc{\calc}{{\cal C}}
\nc{\cald}{\mathcal{D}}
\nc{\cale}{{\cal E}}
\nc{\calf}{{\cal F}}
\nc{\calg}{{\cal G}}
\nc{\calh}{{\cal H}}
\nc{\cali}{{\mathcal I}}
\nc{\call}{{\cal L}}
\nc{\calm}{{\cal M}}
\nc{\caln}{{\cal N}}
\nc{\calo}{{\cal O}}
\nc{\calp}{{\cal P}}
\nc{\calr}{{\cal R}}
\nc{\calt}{{\cal T}}
\nc{\calw}{{\cal W}}
\nc{\calx}{{\mathcal X}}
\nc{\CA}{\mathcal{A}}

\nc{\fraka}{{\frak a}}
\nc{\frakB}{{\frak B}}
\nc{\frakm}{{\frak m}}
\nc{\frakp}{{\frak p}}


\font\cyr=wncyr10

\title{ 
Baxter Algebras, Stirling Numbers and Partitions}
\author{Li Guo} 
\address{
Department of Mathematics and Computer Science,
Rutgers University,\\
Newark, New Jersey 07102, USA}
\email{liguo@newark.rutgers.edu}



\begin{abstract}
Recent developments of Baxter algebras have lead to applications to combinatorics, 
number theory and mathematical physics. 
We relate Baxter algebras to Stirling numbers of the first kind and 
the second kind, partitions and multinomial coefficients. 
This allows us to apply congruences from number theory to obtain congruences 
in Baxter algebras. 
\end{abstract}


\maketitle


\section{Introduction}
\label{intro}

Baxter algebra was first introduced by Baxter~\cite{Ba} in 1960 to study the theory 
of fluctuations in probability. 
Let $C$ be a commutative ring and let $\lambda\in C$ be fixed. 
A Baxter $C$-algebra of weight $\lambda$ is a pair $(R,P)$ in which 
$R$ is a $C$-algebra and $P:R\to R$ is a $C$-linear map such that 
\begin{equation}
 P(x)P(y) =P(xP(y))+P(P(x)y)+\lambda P(xy),\ \forall x,\ y\in R.
\mlabel{eq:Ba}
\end{equation}
The study of Baxter algebras was continued through 1960s and 1970s~\cite{Ro1,Ro2,Ca} 
and has experienced a quite 
remarkable renaissance in the last few years with applications 
to renormalization theory in quantum field theory~\cite{C-K1,C-K2}, 
multiple zeta values in number theory~\cite{Ho} and umbral calculus~\cite{Gu2} 
in combinatorics, as well as in dendriform algebras~\cite{A-L,EF} and  
Hopf algebras~\cite{A-G-K-O}. 
We will give connections of Baxter algebras to Stirling numbers, multinomial 
coefficients and partitions. 


Stirling numbers of the first kind and the second kind were introduced by J. Stirling 
\cite{St} in 1730 and have been studied in the past by well-known mathematicians 
like Euler, Lagrange, Laplace and Cauchy. These numbers play important roles in 
combinatorics, number theory, probability and Statistics. There is a large 
literature on these numbers, as can been seen in the survey article \cite{C-S}.

For each positive integer $n$, let $(t)_n=t(t-1)\cdots(t-n+1)\in \QQ[t]$ be the 
descending (falling) factorial. Also define $(t)_0=1$. Stirling numbers of first 
kind, denoted by $s(n,k)$,  
and Stirling numbers of the second kind, denoted $S(n,k)$ with $n,\ k\in \NN$, 
are defined to be the coefficients in the 
expression 
$$ (t)_n=\sum_{k=0}^n s(n,k) t^k$$
and in the expression 
$$ t^n=\sum_{k=0}^n S(n,k) (t)_k.$$
We will be mostly concerned with Stirling numbers of the second kind. So we 
just mention that the two groups of numbers have similar properties and 
are related by the duality
\begin{equation} \sum_{k=0}^\infty s(n,k)S(k,m)=\sum_{k=0}^\infty S(n,k)s(k,m) =
\delta_{n,m},\ n,\ m\in \ZZ,
\mlabel{eq:dual}
\end{equation}
where $\delta_{n,m}$ is the Kronecker delta. 

Stirling numbers of the second kind can be recursively defined by
\begin{equation}
S(n+1,k+1)=S(n,k)+(k+1)S(n,k+1), \ n,\ k\geq 0
\mlabel{eq:rec}
\end{equation} 
with $S(0,0)=1$, $S(n,0)=S(0,k)=0$ for $n,\ k\geq 1$. 
It follows that $$S(n,1)=S(n,n)=1$$ for $n\geq 1.$ 

Stirling numbers of the second kind also have the following equivalent descriptions: 
\begin{enumerate}
\item 
$S(n,k) = \frac{1}{k!} [\Delta^k t^n]_{t=0},$ where $\Delta$ is the forward 
difference operator.
\item 
$S(n,k)$ is the number of partitions of $n$ objects 
into $k$ non-empty cells. 
\item 
$k! S(n,k)$ is the number of surjective maps from 
a set with $n$ elements to a set with $k$ elements. 
\item 
$\displaystyle{\exp[t(e^u-1)]=\sum_{n=0}^\infty \sum_{k=0}^n 
    S(n,k) t^k u^n/n!.}$
\end{enumerate}

Our first goal is to give an interpretation of Stirling numbers of the first and 
the second kind in terms of Baxter algebras. This is given in Section~\ref{sec:st} after 
introducing the concepts of free Baxter algebras and mixable shuffle products.
In order to ease notations and focus on such connections, we assume 
that the weight $\lambda$ of the Baxter algebra is $1$ in this introduction. 
We will consider Baxter algebras of any weight $\lambda$ in later sections.


\begin{theorem} {\bf $($=Theorem~\ref{thm:conn}$)$}
For any Baxter algebra $(R,P)$ of weight 1 and integer $n\geq 1$, we have 
\begin{equation}
P(1)^n = \sum_{k=1}^n k! S(n,k) P^k(1),\ \ \ \ 
n!P^n(1)=\sum_{k=1}^n s(n,k)  P(1)^k.  
\mlabel{eq:conn}
\end{equation}
\end{theorem}
Recall that $k! S(n,k)$ is the number of ways to put $n$ different objects 
into $m$ different non-empty cells. It is also the number of surjective maps 
from $\{1,\ldots,n\}$ to $\{1,\ldots,k\}$. 
Further, $k! S(n,k)= \Delta^kt^n|_{t=0}$. 

{}From the theorem, we formally have 
$$e^{P(1)u}=\sum_{n=0}^\infty P(1)^n u^n/n! 
    =\sum_{n=0}^\infty \sum_{k=0}^n S(n,k)k! P^k(1) u^n/n!\in R[[u]].$$
So $e^{P(1)u}$ serves as a generating function of $k!S(n,k)$ (compare to item (4) 
in the list before the theorem).
We likewise obtain a generating function for multinomial coefficients 
(see Theorem~\ref{pp:ref}).

Stirling numbers of the second kind enjoy many congruences~\cite{Ca1,Ca2,Co}. 
The most well-known is the congruence 
\begin{equation} S(p,k)\equiv 0 \pmod p,\ 0<k<p 
\mlabel{eq:scong}
\end{equation}
for any prime number $p$. 
There are several proofs of this property~\cite{Co} and new proofs can be found 
as recently as 1997~\cite{H-M-S}. 
As a consequence of (\ref{eq:conn}) and (\ref{eq:scong}), we obtain, for 
any Baxter algebra of weight $\lambda=1$,  
$$P(1)^p \equiv P(1) \mod p.$$

We extend the relation (\ref{eq:conn}) further in Section~\ref{sec:cong} 
and show that multinomial coefficients $\binc{n}{n_1,\cdots,n_k}$ 
and partitions also arise naturally in the products of Baxter  
algebras (Theorem~\ref{pp:ref}). 
This allows us to recover the 
congruence~(\ref{eq:scong}) using multinomial coefficients 
(Proposition~\ref{pp:cong}). 

The congruence (\ref{eq:scong}) is equivalent to Fermat's little theorem 
$$ a^p\equiv a \pmod p,\ a\in \ZZ.$$ 
We apply (\ref{eq:conn}) and (\ref{eq:scong}) in Section~\ref{sec:fermat} 
to generalize this theorem for 
elements in Baxter algebras and obtain 

\begin{theorem} {\bf $($=Theorem~\ref{thm:fermat}$)$}
Let $A$ be a commutative $C$-algebra and 
let $\sha_C(A)$ be the free Baxter algebra on $A$ of weight 1. 
If $a^p\equiv a \mod p$ 
for all $a\in A$, then $a^p\equiv a \mod p$ for all $a\in \sha_C(A)$. 
In particular, taking $C=\FF_p$, 
then $a^p=a$ for all $a\in \sha_{\FF_p}(\FF_p)$.
\end{theorem}

\section{Stirling numbers in Baxter algebras}
\mlabel{sec:st}

We start by reviewing basics on Baxter algebras and the construction of mixed shuffle 
Baxter algebras as a model of free Baxter algebras. 
For further information about Baxter 
algebras, see \cite{G-K,Gu1} 
and references given there. 

Effective study of Baxter algebras depends on explicit constructions of 
free Baxter algebras. This is in analog to the study of commutative algebras 
benefiting from explicit constructions of the free objects (the polynomial algebras). 
For free Baxter algebras we have mixable shuffle algebras. 
We will need the most general construction of mixable shuffle Baxter algebras. 
So we will summarize this construction first. More explicit constructions 
of special cases will be introduced later as needed. 

Given a commutative algebra $C$, $\lambda\in C$ as above, and a $C$-algebra $A$,
the free Baxter $C$-algebra on $A$ is defined to be a Baxter $C$-algebra 
$(\sha_{C,\lambda}(A),P_A)$ together with a $C$-algebra homomorphism 
$j_A:A\to \sha_{C,\lambda}(A)$ with the property that, for any Baxter $C$-algebra 
$(R,P)$ and any $C$-algebra homomorphism $f:A\to R$, there is a unique Baxter 
$C$-algebra homomorphism $\tilde{f}:(\sha_{C,\lambda}(A),P_A)\to (R,P)$ such that 
$j_A\circ \tilde{f}=f$ as $C$-algebra homomorphisms. 

One realization of this free Baxter algebra is given by the mixable shuffle Baxter 
algebra. 
The mixed shuffle Baxter algebra is a pair $(\sha_{C,\lambda}(A),P_A)$ (by misuse 
of notations), where 
$\sha_{C,\lambda}(A)$ is a $C$-algebra in which 
\begin{itemize}
\item
the $C$-module structure is given by 
the direct sum 
$$\bigoplus_{n=1}^\infty A^{\otimes n},\ 
A^{\otimes n}=\underbrace{A\otimes_C\ldots \otimes_C A}_{n-{\rm factors}};$$
\item
the multiplication is given by the mixable shuffle product $\diamond$, recursively defined on 
$A^{\otimes m} \times A^{\otimes n}$ by 
\begin{eqnarray*} 
a_0 \diamond (b_0\otimes b_1 \otimes \ldots \otimes b_n)& =& a_0b_0 \otimes b_1\otimes \ldots 
\otimes b_n, \\
(a_0\otimes a_1\otimes \ldots \otimes a_m)\diamond b_0
& =& a_0b_0 \otimes a_1 \otimes \ldots \otimes a_m, 
\ a_i,\, b_j\in A^{\otimes 1}=A,
\end{eqnarray*}
 and 
\begin{eqnarray}
\lefteqn{(a_0\otimes a_1\otimes \ldots \otimes a_m)\diamond (b_0\otimes b_1 \otimes \ldots 
\otimes b_n) } \notag \\
&=& (a_0 b_0) \otimes \big ((a_1\otimes \ldots \otimes a_m)\diamond 
(1\otimes b_1\otimes \ldots \otimes b_n)\big ) \notag \\
&& +
(a_0 b_0) \otimes \big ((1\otimes a_1 \otimes \ldots \otimes a_m)\diamond (b_1\otimes \ldots 
\otimes b_n)\big) \notag \\
&& + \lambda
a_0 b_0 \otimes \big ((a_1\otimes \ldots \otimes a_m)\diamond (b_1\otimes \ldots \otimes 
b_n)\big),\ a_i,\ b_j\in A.
\mlabel{eq:shuf}
\end{eqnarray}
\end{itemize}
The Baxter operator $P_A$ is defined by 
$$P_A(a_1\otimes \ldots \otimes a_m)=1\otimes a_1\otimes \ldots \otimes a_m,\ 
a_1\otimes \ldots \otimes a_m \in A^{\otimes m},\ 
m\geq 1.$$ 
Since the mixable shuffle product is compatible with the product on $A$, we will 
suppress the notation $\diamond$. We will also express $C$ and $\lambda$ from 
$\sha_{C,\lambda}(A)$ when there is not danger of confusion. 

For example, we have 
\begin{eqnarray*}
\lefteqn{(a_0\otimes a_1\otimes a_2) (b_0\otimes b_1) }\\
&=&a_0b_0\otimes \left ( (a_1\otimes a_2)(1\otimes b_1) 
+ b_1 (1\otimes a_1\otimes a_2)
+ \lambda (a_1\otimes a_2) b_1 \right ) \\
&=& a_0b_0\otimes \left ( a_1 \otimes (a_2 (1\otimes b_1)+(b_1 (1\otimes a_2)) 
    + \lambda a_2 b_1) 
+ b_1 \otimes a_1\otimes a_2
+ a_1 b_1 \otimes a_2 \right ) \\
&=& a_0b_0 \otimes \left ( a_1 \otimes a_2 \otimes b_1 + a_1 \otimes b_1 \otimes a_2 
    + \lambda a_1 \otimes a_2 b_1 + b_1 \otimes a_1 \otimes a_2 + \lambda a_1 b_1 \otimes a_2 
    \right ). 
\end{eqnarray*}



\begin{theorem}
Let $(R,P)$ be a Baxter $C$-algebra. For $n\geq 1$, 
\begin{equation}
 P(1)^n=\sum_{k=1}^n k! S(n,k) \lambda^{n-k} P^k(1),\ \ \ \ 
n!P^n(1)=\sum_{k=1}^n s(n,k) \lambda^{n-k} P(1)^k.  
\mlabel{eq:con}
\end{equation}
\mlabel{thm:conn}
\end{theorem}

\proofbegin 
The theorem can be proved using properties of free Baxter algebras. 
We will give a more direct proof by induction. 
For the first equation, it obviously holds for $n=1$ 
since $S(1,1)=1$. 
Assume that it is true for $n$. 
By~\cite[Prop. 6.1]{G-K}, 
\begin{equation}
P^k(1) P(1)=(k+1)P^{k+1}(1) + k \lambda P^{ k}(1).
\mlabel{eq:cir}
\end{equation}
Then we have 

\begin{eqnarray*}
\lefteqn{ P(1)^{n+1}= P(1)^n P(1)}\\ 
&=& \big (\sum_{k=1}^n k! S(n,k) \lambda^{n-k} P^{ k}(1) \big )P(1) \\ 
&=& \sum_{k=1}^n k! S(n,k) \lambda^{n-k} 
    \big ((k+1)P^{ k+1}(1) +\lambda k P^{ k}(1)\big ) \\
&=& \sum_{k=1}^n (k+1)! S(n,k) \lambda^{n-k} 
    P^{ k+1}(1) 
    +\sum_{k=1}^n k! S(n,k) \lambda^{n+1-k} 
    k P^{ k}(1). 
\end{eqnarray*}
Replacing $k$ by $k-1$ in the first sum 
we have 
\begin{eqnarray*}
P(1)^{n+1}&=&\sum_{k=2}^{n+1} k! S(n,k-1) 
    \lambda^{k+1-n} P^{ k} 
    +\sum_{k=1}^n k! S(n,k) \lambda^{k+1-n} 
    k P^{ k}\\  
&=& S(n,1) \lambda^k P(1)
    +\sum_{k=2}^n k! \lambda^{n+1-k} 
    \big (S(n,k-1) +k S(n,k) P^{ k}(1)\big )\\
&&    +(n+1)! S(n,n) P^{ n+1}(1). 
\end{eqnarray*}
So using $S(k,1)=S(k,k)=1,\ k\in\NN$ and (\ref{eq:rec}), we obtain
\begin{eqnarray*}
P(1)^{n+1}&=& S(n+1,1) \lambda^n P(1)
    +\sum_{k=2}^n k! \lambda^{n+1-k} 
    S(n+1,k) P^{ k} (1) \\
&&    +(n+1)! S(n+1,n+1) P^{ (n+1)}(1) \\
&=&   \sum_{k=1}^{n+1} k! \lambda^{n+1-k} 
    S(n+1,k) P^{ k}(1).  
\end{eqnarray*}
This proves the first equation. The second equation then follows from 
the duality in Eq.~(\ref{eq:dual}) between Stirling numbers of the first and the second kind.
\proofend

It follows from Theorem~\ref{thm:conn} and the congruence~(\ref{eq:scong}) 
that 
\begin{coro}
Let $(R,P)$ be a Baxter $C$-algebra and let $p$ be a prime number. 
\begin{equation}
P(1)^p \equiv \lambda^{p-1} P(1) \pmod p
\mlabel{eq:bcong1}
\end{equation}
\end{coro}
A natural question to ask is whether (\ref{eq:bcong1}) holds when 
$1$ is replaced by other elements of $R$. This question will be addressed 
in Theorem~\ref{thm:fermat}.

\section{Multinomial coefficients and partitions in Baxter algebras}
\mlabel{sec:cong}
Let $X=\{x\}$. Let 
$$\cali=\{I=(i_1,\cdots,i_k)| i_s\geq 0, 1\leq s\leq k, k\geq 1\}.$$
For $I\in \cali$, define the norm of $I$ to be 
$|I|=\sum_{s=1}^k i_s$ and define the length of $I$ to be 
$\ell(I)=k$. 
Define $I>0$ if $i_s>0$ for $1\leq s\leq k$. 

Denote $$x^{\otimes I}=x^{i_1}\otimes \cdots \otimes x^{i_k}$$ 
and 
$$\calx=\big\{x^{\otimes I}\big |I\in \cali\big\}.$$ 
Recall that the free Baxter algebra on $C[x]$ is 
$\bigoplus_{n\geq 1} C[x]^{\otimes n}$ as a $C$-module, 
so is a free $C$-module on $\calx$ together 
with the mixed shuffle product defined in (\ref{eq:shuf})~\cite{Gu1,G-K}. 
The Baxter operator $P_X$ on $\sha_C(X)$ is defined by $P_X(x^{\otimes I})=
1\otimes x^{\otimes I}.$ 
Thus $P_X(\sha_C(X))=1\otimes \sha_C(X)$ is a free $C$-module on 
$1\otimes \calx:=\big \{ 1\otimes x^{\otimes I}\ \big |\ I\in \cali\big \}.$

By the Baxter identity~(\ref{eq:Ba}), $P_X(\sha_C(X))$ is closed under multiplication. 
Therefore, for each $n\geq 1$, 
\begin{equation}
(1\otimes x)^n=\sum_{|I|=n} S(n,I)\otimes x^{\otimes I}
\mlabel{eq:refine}
\end{equation}
for elements $S(n,I)\in C$. 
The first few terms of $S(n,I)$ are given by 
\begin{eqnarray*}
(1\otimes x)&=&1\otimes x,\\
(1\otimes x)^2 &=& 2 (1\otimes x\otimes x)+\lambda (1\otimes x^2),\\
(1\otimes x)^3 &=& 6 (1\otimes x\otimes x\otimes x) + 3\lambda (1\otimes x\otimes x^2) 
    +3\lambda (1\otimes x^2\otimes x) + \lambda^2(1\otimes x^3),\\
(1\otimes x)^4 &=& 24 (1\otimes x\otimes x\otimes x\otimes x)+
    12\lambda (1\otimes x\otimes x \otimes x^2) 
    +12\lambda (1\otimes x\otimes x^2\otimes x)\\
&&    +12\lambda (1\otimes x^2\otimes x\otimes x)+6\lambda^2 (1\otimes x^2\otimes x^2) 
    +4\lambda^2 (1\otimes x \otimes x^3)\\
&&  +4\lambda^2 (1\otimes x^3\otimes x)+ \lambda^3(1\otimes x^4).
\end{eqnarray*}
\begin{lemma}
$S(n,I)\neq 0$ only if $I>0$. 
\mlabel{lem:pos}
\end{lemma}
\begin{proof}
We prove by induction on $n$ with $n=1$ being obvious. 
Assuming the lemma holds for $n$. So 
$$(1\otimes x)^n=\sum_{|I|=n, I>0} S(n,I)\otimes x^{\otimes I}.$$
For $x^{\otimes I}=x^{i_1}\otimes\cdots\otimes x^{i_k}$, the 
mixable shuffle product (\ref{eq:shuf}) gives
\begin{eqnarray*}
(1\otimes x^{\otimes I})(1\otimes x)&=& 
1\otimes \big (x^{(i_1,\cdots,i_k,1)}+x^{\otimes (i_1,\cdots,i_{k-1},1,i_k)} 
+\cdots + x^{\otimes (1,i_1,\cdots,i_k)} \\
&&+\lambda x^{\otimes (i_1,\cdots,i_{k-1},i_k+1)}
+\cdots + \lambda x^{\otimes (i_1+1,i_2,\cdots,i_k)}\big ).
\end{eqnarray*}
This shows that the lemma holds for $n+1$. 
\end{proof}

By Lemma~\ref{lem:pos}, the sum in (\ref{eq:refine}) is over the set of {\em ordered} 
partitions of $n$.  For example, when $n=3$, the terms on the right of 
$(1\otimes x)^3$ correspond to the ordered partitions
$$n=1+1+1=1+2=2+1=3.$$
For the coefficients in the sum, we have 

\begin{theorem} 
\begin{enumerate}
\item For $I>0$, we have $\displaystyle{S(n,I)=\binc{n}{i_1,\cdots,i_k}} \lambda^{n-\ell(I)}$. 
\item $\displaystyle{\sum_{|I|=n,\ell(I)=k} S(n,I)}=k!S(n,k)\lambda^{n-k}.$
\end{enumerate}
\mlabel{pp:ref}
\end{theorem}
Recall that the multinomial coefficient $\binc{n}{i_1,\cdots,i_k}$ is 
the number of ways to distribute $n$ distinct objects into $k$ distinct 
sets $S_1,\cdots,S_k$ with $i_1,\cdots,i_k$ objects. 

By Eq. (\ref{eq:refine}) and the above theorem, we formally have 
$$ \frac{1}{1-(1\otimes x)} =\sum_{n=0}^\infty (1\otimes x)^n 
    =\sum_{n=0}^\infty \sum_{|I|=n} \binc{n}{I} \otimes x^I.$$
Here $\binc{n}{I}=\binc{n}{i_1,\cdots,i_k}$ for $I=(i_1,\cdots,i_k)$ by abbreviation 
and $(1\otimes x)^0=1$ by convention. Thus $1/(1-1\otimes x)$ can be regarded as 
a generating function for the multinomial coefficients $\binc{n}{I}$.

\proofbegin (1). 
We prove by induction on $\ell(I$. It is obviously true for $\ell(I)=1$. 
Assume that it is true for all $\ell(I)$ with $\ell(I)=n$. 
Then 
\begin{eqnarray*}
\sum_{|J|=n+1} S(n+1,J)\lambda^{n+1-\ell(J)} (1\otimes x^{\otimes J})
&=& (1\otimes x)^{n+1}\\
&=& \sum_{|I|=n} S(n,I)(1\otimes x^{\otimes I})(1\otimes x).
\end{eqnarray*}

Let $I=(i_1,\cdots,i_k)$. Using the notation 
$x^{\otimes I}=x^{i_1}\otimes\cdots\otimes x^{i_k}$, the 
mixable shuffle product (\ref{eq:shuf}) gives
\begin{eqnarray*}
(1\otimes x^{\otimes I})(1\otimes x)&=& 
1\otimes \big (x^{(i_1,\cdots,i_k,1)}+x^{\otimes (i_1,\cdots,i_{k-1},1,i_k)} 
+\cdots + x^{\otimes (1,i_1,\cdots,i_k)} \\
&&+\lambda x^{\otimes (i_1,\cdots,i_{k-1},i_k+1)}
+\cdots + \lambda x^{\otimes (i_1+1,i_2,\cdots,i_k)}\big ).
\end{eqnarray*}
It follows  that,  
the coefficient of $J=(j_1,\cdots,j_k)$ with $|J|=j_1+\cdots+j_k=n+1$ 
has contribution exactly from $I_t$ defined by 
$I_t=(j_1,\cdots,j_{t-1},j_t-1,j_{t+1},\cdots,j_k),\ 1\leq t\leq k,$
with the convention that $i_t=j_t-1$ is deleted from $I_t$ 
in the event that $j_t=1$. 
Then $\ell(I_t)=\ell(J)$ if $j_t>1$ and $\ell(I_t)=\ell(J)-1$ if $j_t=1$. 
For example, when $J=(2,1,4)$, we have $I_1=(1,1,4),\ I_2=(2,4),$ and
$I_3=(2,1,3).$
We then have 
\begin{eqnarray*}
S(n+1,J)\lambda^{n+1-\ell(J)}&=& \sum_{j_t=1} S(n,I_t)\lambda^{n-\ell(I_t)}
    +\sum_{j_t>1} S(n,I_t) \lambda^{n-\ell(I_t)+1} \\
&=& \sum_{j_t=1} S(n,I_t)\lambda^{n+1-\ell(J)}
    +\sum_{j_t>1} S(n,I_t) \lambda^{n+1-\ell(J)} \\
&=& \sum_{t=1}^k S(n,I_t)\lambda^{n+1-\ell(J)}. 
\end{eqnarray*}
So by induction we have 
\begin{eqnarray*}
S(n+1,J)\lambda^{n+1-\ell(J)}
&=& \sum_{t=1}^k \binc{n}{j_1,\cdots,j_{t-1},j_t-1,j_{t+1},\cdots,j_k}\\
&=& \binc{n+1}{j_1,\cdots,j_k}.
\end{eqnarray*}
Where the last equation is the generalized Pascal identity~\cite[p.24]{Ri}. 

(2). 
Once part 1 is proved, part 2 is a standard property on Stirling numbers and 
multinomial coefficients. We give a direct proof using properties of free 
Baxter algebras instead of part 1. 
Since $\sha(\{x\})$ is the free Baxter algebra on $\{x\}$, there is a unique 
Baxter algebra homomorphism 
$$f: \sha(\{x\}) \to \sha(C)$$ induced by $f(x)=1$. 
We have 
\begin{eqnarray*}
(1\otimes 1)^n&=& f((1\otimes x)^n)\\
&=& \sum_{|I|=n} f(S(n,I)\otimes x^{\otimes I})\\
&=& \sum_{|I|=n} S(n,I)\otimes f(x^{\otimes I})\\
&=& \sum_{|I|=n} S(n,I)\otimes f(x)^{\otimes I}\\
&=& \sum_{|I|=n} S(n,I)\otimes 1^{\otimes I}\\
&=& \sum_{|I|=n} S(n,I)\otimes 1^{\otimes \ell(I)}\\
&=& \sum_{|I|=n} \sum_{\ell(I)=k} S(n,I)\otimes 1^{\otimes k}.
\end{eqnarray*}
On the other hand, by Theorem~\ref{thm:conn}, 
$$(1\otimes 1)^n=P_X(1)^n=\sum_{k=1}^n k! S(n,k) \lambda^{n-k} (1\otimes 1)^k.$$
Since $\{1\otimes 1^{\otimes k}\big |\ k\geq 1\}$  is linearly independent 
over $C$, we obtain 
$$\sum_{|I|=n, \ell(I)=k} S(n,I)=k!\,S(n,k)\lambda^{n-k},$$ 
as needed. 
\proofend

We now deduce congruences of $S(n,I)$. As a result, 
we obtain a proof of the classical congruence (\ref{eq:scong}). 

\begin{prop} Let $p$ be a prime number. 
\begin{enumerate}
\item $p| S(p,I)$ for any $I\in \cali$ with $|I|=p$ and $I\neq (p)$. 
\item $p| S(p,k)$ for $1<k<p$. 
\item $(1\otimes x)^p \equiv \lambda^{p-1}\otimes x^p \pmod{p}.$
\item For any Baxter algebra $(R,P)$ and $a\in R$, we have 
$P(a)^p\equiv \lambda^{p-1} P(a^p) \pmod{p}.$
\item Let $\lambda^{p-1}=\lambda$ and let $\sha_{\FF_p}(\FF_p)$ be the 
free Baxter algebra over $\FF_p$ of weight $\lambda$.  
For $a\in \sha_{\FF_p}(\FF_p)$, we have $a^p=\lambda^{p-1}a$. 
\end{enumerate}
\mlabel{pp:cong}
\end{prop}
Part 5 of Proposition~\ref{pp:cong} can be regarded as the Fermat's little 
theorem for Baxter algebra in a special case. We will consider the general 
case in the next section. 

\proofbegin
(1). 
By Theorem~\ref{pp:ref},
 $S(p,I)=\disp{\frac{p!}{n_1!\cdots n_k!}}\lambda^{p-k}$. For $I\neq (p)$, 
we have $p\nmid n_1!\cdots n_k!$ since $n_i<p$. 
The claim follows. 

\noindent
(2). From (1) and Theorem~\ref{pp:ref} (2), we have 
$p|k!S(p,k)\lambda^{p-k},\ 1<k<p$. So we just need to take $\lambda=1$. 

\noindent
(3). This follows from
$$ (1\otimes x)^p = \sum_{|I|=p,\ell(I)=k} S(p,I)\lambda^{p-\ell(k)}\otimes x^{\otimes I} 
\equiv \lambda^{p-1}\otimes x^{\otimes (p)}=\lambda^{p-1}\otimes x^p \mod p.$$

\noindent
(4). Since $\sha(\{x\})$ is the free Baxter algebra on $\{x\}$, 
the assignment $x\mapsto a$ induces a Baxter algebra homomorphism 
$f:\sha(\{x\}) \to R.$ Applying $f$ to the congruence in (3) gives (4). 

\noindent 
(5). We first prove $(1^{\otimes n})^p=1^{\otimes n}$ for $n\geq 1$ 
with the case of $n=1$ being clear.  
Using (4) and induction, we have 
\begin{eqnarray*}
(1^{\otimes (n+1)})^p&=& P(1^{\otimes n})^p\\
&=& \lambda^{p-1} P((1^{\otimes n})^p)\\
&=& P((1^{\otimes n})^p)\\
&=& P(1^{\otimes n})\\
&=& 1^{\otimes (n+1)}.
\end{eqnarray*}
Since $\{1^{\otimes n} \big | n\geq 1\}$ is a basis 
of $\sha_{\FF_p}(\FF_p)$ and $\sha_{\FF_p}(\FF_p)$ has characteristic $p$, 
(5) follows. 
\proofend

\section{Fermat's little theorem for Baxter algebras}
\mlabel{sec:fermat}
Now we consider free Baxter algebras in general. 
Let $A$ be $C$-algebra. 
Recall that the free Baxter algebra $\sha_C(A)$ is 
$$\sha_C(A)=\bigoplus_{n\geq 1} A^{\otimes n}$$ 
with multiplication given by the mixable shuffle product and with Baxter operator
$P_A(a_1\otimes \cdots \otimes a_n)=1\otimes a_1\otimes \cdots \otimes a_n$ 
for $a_1\otimes \cdots \otimes a_n\in A^{\otimes n}.$ An Baxter ideal 
of $\sha_C(A)$ is defined to be an ideal $I$ of $\sha_C(A)$ such that 
$P(I)\subseteq I$.

\begin{theorem} Let $p$ be a prime number. Assume that the weight $\lambda$ 
of the Baxter algebra is invertible. 
\begin{enumerate}
\item 
For $a_1\otimes \cdots \otimes a_n\in A^{\otimes n}$, 
$$(a_1\otimes \cdots \otimes a_n)^p 
    \equiv \lambda^{(n-1)(p-1)} a_1^p\otimes \cdots \otimes a_n^p 
    \mod p.$$
In particular, if $\lambda=1$, then 
$$(a_1\otimes \cdots \otimes a_n)^p 
    \equiv a_1^p\otimes \cdots \otimes a_n^p 
    \mod p.$$
\item 
Let $I$ be the Baxter ideal of $\sha_C(A)$ generated by $p$ and 
$\{a^p-a \big |\ a\in A\}$. Then for $a\in \sha_C(A)$, $a^p \equiv a \mod I.$
\item
Let $A$ be a $C$-algebra of characteristic $p$ such that $a^p=a$ for all $a\in A$. 
Then $a^p=a$ for all $a\in \sha_C(A)$. 
\end{enumerate}
\mlabel{thm:fermat}
\end{theorem}

\begin{remark}
\begin{enumerate}
\item Taking a tensor product and taking a power cannot be interchanged in general 
in $\sha_C(A)$. This is similar to taking a sum and taking a power. So 
(1) can be regarded as a tensor form of the freshman's dream 
$(a+b)^p\equiv a^p+b^p\mod p.$
\item
Let $A$ be the $\FF_p$-algebra $\FF_p[x]/(x^p-x)$, then 
$a^p=a$ for all $a\in A$. This follows from Fermat's little theorem and 
the freshman's dream. 
The free Baxter algebra $\sha_{\FF_p}(A)$ is a quite large algebra over 
$A$, 
$$\sha_{\FF_p}(A)=\bigoplus_{n\geq 1} A^{\otimes n}.$$
But we still have 
$a^p= a$ for all $a\in \sha_{\FF_p}(A)$ by the theorem.  
\end{enumerate}
\end{remark}

\proofbegin
(1). By (4) of Proposition~\ref{pp:cong}, we have 
\begin{eqnarray*}
(a_1\otimes \cdots  \otimes a_n)^p &\equiv & 
(a_1P_A(a_2\otimes \cdots \otimes a_n))^p \\
&\equiv & a_1^p P_A(a_2\otimes \cdots \otimes a_n)^p \\
&\equiv & a_1^p \lambda^{p-1} P_A((a_2\otimes \cdots \otimes a_n)^p) \\
&\equiv & \lambda^{p-1} a_1^p\otimes (a_2\otimes \cdots \otimes a_n)^p\mod p.
\end{eqnarray*}
So by induction, we have 
$$ (a_1\otimes \cdots  \otimes a_n)^p\equiv \lambda^{(p-1)(n-1)} a_1^p\otimes a_2^p
\otimes \cdots \otimes a_n^p \mod I.$$

\noindent
(2). By part 1 we have
$$(a_1\otimes \cdots \otimes a_n)^p \equiv 
\lambda^{(n-1)(p-1)} a_1^p \otimes \cdots \otimes a_n^p
\equiv \lambda^{(n-1)(p-1)} a_1 \otimes \cdots \otimes a_n \mod I.$$
Since $\lambda$ is invertible, from $\lambda^p\equiv \lambda \mod I$ we have 
$\lambda^{p-1}\equiv 1\mod I$. (2) follows. 

\noindent
(3). This follows from (2), since the ideal $I$ is trivial here. 
\proofend

\section*{Acknowledgements}
The author thanks Andrew Granville for helpful discussions on Stirling numbers.

\end{document}